# METHOD OF TRAINING EXAMPLES IN SOLVING INVERSE ILL-POSED PROBLEMS OF SPECTROSCOPY

## V.S. Sizikov, A.V. Stepanov

Further development of the method of computational experiments for solving ill-posed problems is given. The effective (unoverstated) estimate for solution error of the first-kind equation is obtained using the truncating singular numbers spectrum of an operator. It is proposed to estimate the magnitude of the truncation by results of solving model (training, learning) examples close to the initial example (problem). This method takes into account an additional information about the solution and gives a new principle for choosing the regularization parameter and error estimate for equation solution by the Tikhonov regularization method. The method is illustrated by a numerical example from the inverse problem of spectroscopy.

**Keywords:** ill-posed problems, Tikhonov regularization, solution error, method of training examples, inverse problem of spectroscopy, integral equation, spread function of spectral device, measured spectrum, training spectra, restored spectrum.

## Introduction

In this paper, we develop an adaptive method of computational experiments for estimating the solution error and choosing the regularization parameter in solving ill-posed problems by the Tikhonov regularization method. Other names: the technique of model, standard, learning, training examples, the way of the pseudoinverse operator [1–8]. This method takes into account an additional (a priori) information about the desired solution (an estimate of the number of maxima, their abscissas and ordinates, etc.) and, in this respect, resembles the methods such as the Tikhonov $\alpha$-regularization with constraints on the solution [9], solution on a compact [4, 9], the methods of descriptive regularization [10], also taking into account a priori information on the solution (nonnegativity, monotonicity, convexity, parameters of extrema, etc.). However, the specific implementation of the method of computational experiments differs from these methods.

This method was earlier developed and applied to signal processing [1–3], image restoration [5, 6] and spectroscopy [5–8]. In this paper, we propose its modification and the application to the inverse problem of spectroscopy.

## Basic relations

Consider an operator equation of the first kind

$$Ay = f, \quad y \in Y, \quad f \in F, \tag{1}$$

where $y$ is desired, and $f$ is given elements of Hilbert spaces $Y$ and $F$; $A$ is a linear bounded operator from $Y$ into $F$. The operator $A$ is not expected to be continuously invertible, i.e. the problem of solving (1) is ill-posed. However, for the exact $f$ we assume that equation (1) is solvable.

The problem is to find an element $y \in Y$ with minimal norm, which supplies the minimum value for the discrepancy $\|Ay - f\|_F$ and which is the pseudosolution, in particular, the normal solution [4, 9, 11].

In the zero-order Tikhonov regularization method [4, 9, 11, 12], giving one of the most effective ways for obtaining pseudosolutions, instead of (1) the equation

$$(\alpha E + \widetilde{B}) y_\alpha = \widetilde{A}^* \widetilde{f}, \tag{2}$$

is solved, where

$$\widetilde{A} = A + \Delta A, \quad \widetilde{f} = f + \Delta f, \quad y_\alpha = y + \Delta y_\alpha, \tag{3}$$

moreover, $A$, $f$ and $y$ are the exact operator and elements; $\widetilde{A}$, $\widetilde{f}$ and $y_\alpha$ are their practical values; $\Delta A$, $\Delta f$ and $\Delta y_\alpha$ are their errors; $\alpha > 0$ is the regularization parameter; $\widetilde{B} = \widetilde{A}^* \widetilde{A}$; $E$ is the unit operator.

**Estimate of solution error.** Consider the question of estimating the error $\Delta y_\alpha$ of the regularized solution $y_\alpha$ and choosing the regularization parameter $\alpha$.

It is known [13] that it is almost impossible to obtain an effective (unoverstated) estimate of the error $\Delta y_\alpha$ without using additional (a priori) information on the solution. In this paper, we propose to use the results of solving "close" model, learning examples as additional information. Taking into account the ratio $\widetilde{A}^* A y = \widetilde{A}^* f$, which follows from (1), as well as the ratios (2) and (3), we obtain

$$(\alpha E + \widetilde{B})\Delta y_\alpha = -(\alpha E + \widetilde{A}^* \Delta A) y + \widetilde{A}^* \Delta f$$

or

$$\Delta y_\alpha = (\alpha E + \widetilde{B})^{-1} \widetilde{A}^* (\Delta f - \Delta A\, y) - \alpha(\alpha E + \widetilde{B})^{-1} y,$$

from where we find the following estimates in the norm of absolute and relative errors of the regularized solution

$$\|\Delta y_\alpha\| \leq \|(\alpha E + \widetilde{B})^{-1} \widetilde{A}^*\| \cdot (\|\Delta f\| + \|\Delta A\| \cdot \|y\|) + \alpha \|(\alpha E + \widetilde{B})^{-1}\| \cdot \|y\|,$$

$$\frac{\|\Delta y_\alpha\|}{\|y\|} \leq \|(\alpha E + \widetilde{B})^{-1} \widetilde{A}^*\| \cdot \left(\frac{\|\Delta f\|}{\|y\|} + \|\Delta A\|\right) + \alpha \|(\alpha E + \widetilde{B})^{-1}\|. \quad (4)$$

Taking into account that $\|A\| \cdot \|y\| \geq \|f\|$ or $1/\|y\| \leq \|A\|/\|f\|$, we obtain the estimate (4) for the relative error of the regularized solution in the form

$$\sigma_{rel}(\alpha) \equiv \frac{\|\Delta y_\alpha\|}{\|y\|} \leq \|(\alpha E + \widetilde{B})^{-1} \widetilde{A}^*\| \cdot \|A\| (\delta_{rel} + \xi_{rel}) + \alpha \|(\alpha E + \widetilde{B})^{-1}\|, \quad (5)$$

where

$$\delta_{rel} = \frac{\|\Delta f\|}{\|f\|}, \quad \xi_{rel} = \frac{\|\Delta A\|}{\|A\|}$$

are the relative errors of the right-hand side $f$ and operator $A$. The right-hand side of (5) is the upper envelope of the true relative error $\sigma_{rel}(\alpha)$. The first summand in the right-hand side of (5) is due to the errors of data, while the second summand is determined by regularization. In (5), we have (cf. [14–17]): $\|(\alpha E + \widetilde{B})^{-1} \widetilde{A}^*\| \leq 1/(2\sqrt{\alpha})$, and norm $\|(\alpha E + \widetilde{B})^{-1}\|$ can be expressed through minimum singular number $\mu$ of symmetric positively determined operator $\alpha E + \widetilde{B}$:

$$\|(\alpha E + \widetilde{B})^{-1}\| = 1/\mu_{min}(\alpha E + \widetilde{B}) = 1/(\alpha + \mu_{min}(\widetilde{B})).$$

We obtain (cf. [14–17]):

$$\sigma_{rel}(\alpha) \equiv \frac{\|\Delta y_\alpha\|}{\|y\|} \leq \frac{\|\widetilde{A}\|}{2\sqrt{\alpha}} \eta + \frac{\alpha}{\alpha + \mu_{min}(\widetilde{B})}, \quad (6)$$

where $\eta = \delta_{rel} + \xi_{rel}$.

However, in practice, estimate (6) (as well as (5)) may give a significant overstatement for $\sigma_{rel}(\alpha)$, since, in the case of ill-conditioned and ill-posed problems, $\mu_{min}(\widetilde{B})$ is close or equal to zero and then (when $\mu_{min}(\widetilde{B}) = 0$)

$$\sigma_{rel}(\alpha) \equiv \frac{\|\Delta y_\alpha\|}{\|y\|} \leq \frac{\|\widetilde{A}\|}{2\sqrt{\alpha}} \eta + 1. \quad (7)$$

The estimate (7) is not only overstated, but not having the minimum with respect to $\alpha$.

To obtain more effective estimate of $\sigma_{rel}(\alpha)$ we use the concept of the pseudoinverse operator having enclosed in it, however, a sense somewhat different from the pseudo-inverse Moore-Penrose matrix $A^+$, giving the solution $y = A^+ f$ [4, 12, 18], and from the regularized operator $(\alpha E + \widetilde{B})^{-1} \widetilde{A}^*$, giving the solution $y_\alpha = (\alpha E + \widetilde{B})^{-1} \widetilde{A}^* f$. The fact is that $A^+$ cor-responds to the case $\alpha \to 0$ and $\mu_{\min}(\widetilde{B}) \neq 0$, while regularization is dealing with a finite value of $\alpha > 0$ and $\mu_{\min}(\widetilde{B}) \approx 0$, which leads to an overstatement of $\sigma_{\rm rel}(\alpha)$ in both cases.

In order to bring the estimate $\sigma_{\rm rel}(\alpha)$ nearer to the true estimate of $\sigma_{\rm rel}(\alpha)$, we *truncate* the spectrum of the operator (matrix in the discrete case) $\widetilde{B}$ from below, namely, instead of $\mu_{\min}(\widetilde{B})$ we use a value $g > \mu_{\min}(\widetilde{B})$ and write (6) in the form

$$\sigma_{\rm rel}(\alpha) \equiv \frac{\|\Delta y_\alpha\|}{\|y\|} \leq \varepsilon(\alpha), \tag{8}$$

where

$$\varepsilon(\alpha) = \frac{\|\widetilde{A}\|}{2\sqrt{\alpha}} \eta + \frac{\alpha}{\alpha + g}. \tag{9}$$

It was shown in [1, 2, 4] that the function $\varepsilon(\alpha)$, according to (9), has a (unique) minimum under the condition

$$\frac{\|\widetilde{A}\|}{\sqrt{g}} \eta < \frac{3\sqrt{3}}{4} \approx 1.30. \tag{10}$$

From the condition $\varepsilon'(\alpha) = 0$, we obtain the equation for $\alpha$ (cf. [2, 4])

$$\alpha = \left(\frac{\|\widetilde{A}\|\eta}{4g}\right)^{2/3} (\alpha + g)^{4/3}. \tag{11}$$

As shown in [2], in this case $\varepsilon''(\alpha) > 0$, i.e. (11) corresponds to the minimum of the function $\varepsilon(\alpha)$.

According to relations (8) and (9), an relative error estimate $\|\Delta y_\alpha\|/\|y\|$ of regularized solution $y_\alpha$ depends on $\widetilde{A}$ and $\eta$ (more exactly, on the product $\|\widetilde{A}\|\eta$). Therefore, if we solve a few examples (e.g., process a few spectra) with the same $\widetilde{A}$ and $\eta$ (spread function and noise), then their error estimates (9) will be identical and unoverstated (in function of $\alpha$). It follows that when solving some original example $P$ (i.e. when processing $\widetilde{f}_P$) with unknown solution (spectrum) $y_P$, one can use the results of solving other (model, training) example $Q$ with known (given) exact solution (spectrum) $y_Q$, with the same $\widetilde{A}$ and $\eta$ as in example $P$. Furthermore, when solving example $Q$, one can calculate the function $\sigma_{\rm rel}(\alpha)_Q = \|\Delta y_{\alpha Q}\|/\|y_Q\|$ and, based on this function, find $\alpha_{{\rm opt}\, Q}$ (optimal value of $\alpha$, at which $\sigma_{\rm rel}(\alpha)_Q = \min_\alpha$). This value $\alpha_{{\rm opt}\, Q}$ can be used for solving the original example (spectrum) $P$.

**Estimate of parameter $g$.** Furthermore, it is necessary to determine the parameter $g$ (see (9)). An estimate of $g$ can be obtained *graphically*, namely, by fitting such value of $g$, at which envelope $\varepsilon(\alpha)$ contacts curve (or a set of curves) $\sigma_{\rm rel}(\alpha)_Q$ (see later Fig. 2). The value of $\alpha$ corresponding to the contact point we denote as $\alpha_g$.

Determining $g$ can also be performed *analytically*. Equating $\sigma_{\rm rel}(\alpha)$ and $\varepsilon(\alpha)$, as well as taking into account the condition $\varepsilon'(\alpha) = 0$, we obtain two equations for two unknowns $\alpha$ and $g$:

$$\left.\begin{aligned}\frac{\|\widetilde{A}\|}{2\sqrt{\alpha}} \eta + \frac{\alpha}{\alpha + g} &= \sigma_{\rm rel}(\alpha), \\ \alpha &= F(\alpha),\end{aligned}\right\} \tag{12}$$

where (see (11))

$$F(\alpha) = \chi(\alpha + g)^{4/3}, \quad \chi = \left(\frac{\|\widetilde{A}\|\eta}{4g}\right)^{2/3}. \tag{13}$$

Here, $\sigma_{rel}(\alpha)$ is the calculated upper curve from a set of curves $\sigma_{rel}(\alpha)_Q = \|\Delta y_{\alpha Q}\| / \|y_Q\|$. The first equation in (12) is the condition of contact of $\varepsilon(\alpha)$ (see (9)) and $\sigma_{rel}(\alpha)$, whereas the second equation is the minimum condition of function $\varepsilon(\alpha)$, i.e. $\varepsilon'(\alpha) = 0$ at the contact point. The first equation can be resolved relatively to $g$:

$$g = \alpha \cdot \left[\left(\sigma_{rel}(\alpha) - \frac{\|\widetilde{A}\|\eta}{2\sqrt{\alpha}}\right)^{-1} - 1\right]. \tag{14}$$

Then, obtained system of two equations can be solved by iterations:

$$\alpha_0 = \chi, \quad \alpha_i = \chi_{i-1}(\alpha_{i-1} + g_{i-1})^{4/3}, \quad \chi_{i-1} = \left(\frac{\|\widetilde{A}\|\eta}{4g_{i-1}}\right)^{2/3},$$

$$g_i = \alpha_i \cdot \left[\left(\sigma_{rel}(\alpha_i) - \frac{\|\widetilde{A}\|\eta}{2\sqrt{\alpha_i}}\right)^{-1} - 1\right], \quad i = 1,2,3\ldots \tag{15}$$

This iterative process for $\alpha$ converges to some $\alpha = \alpha_g$, since $|F'(\alpha)| < 1$.

However, since the function $\sigma_{rel}(\alpha)$ is given in tabular form, it is more convenient to solve the problem *graphically* displaying onto a computer monitor the curves $\sigma_{rel}(\alpha)$ and $\varepsilon(\alpha)$ at different $g$ (Fig. 2). To enhance the efficiency of this method when working out a model example $Q$ (or several examples) it is necessary to use an additional information about the original example (spectrum) $P$, namely, an estimate of the number of maxima (spectral lines) in the desired solution (spectrum) $y_P$, ratios of their intensities, values of its abscissa (wavelengths or frequencies), the type of kernel (SF), etc. Use of such information will allow to choose more "successfully" the regularization parameter $\alpha$ and estimate the error of solving the examples (spectra) $Q$ and $P$.

The modeling method generates a *regularizing algorithm* (RA), since when $\eta \to 0$, $\alpha = o(\eta^2)$ and finite $\|\widetilde{A}\|$ and $g$, we have for original and model examples according to (8) and (9) (cf. [2,4]):

$$\sigma_{rel}(\alpha) \equiv \frac{\|\Delta y_\alpha\|}{\|y\|} \to 0,$$

i.e. at zero errors of initial data, the solution $y_\alpha$ turns into the exact solution (normal pseudo-solution).

**Remark 1.** Although the method of modeling (training, learning) requires a lot of preliminary work on drawing up and solving the model (training) examples, it is very effective in cases when it is required to solve a significant number of "close" examples (to resolve signals for a number of times, to restore several similar spectra in the inverse spectroscopy problem, etc.). Moreover, this method allows for a number of training examples to explore practical possibilities of the used method and algorithm applied to a particular problem (to obtain the real solution error, the possibility of restoring the fine solution structure, etc.).

**Remark 2.** The objection can arise that, because of the ill-posedness of a problem, even small differences of the model problem (example, spectrum) from the original one can lead to significant differences of the regularization parameter $\alpha$, the relative solution error $\|\Delta y_\alpha\| / \|y\|$, etc. However, firstly, the problem is solved by a stable regularization method and it is the conditionally well-posed (by Tikhonov), and secondly, relations (8) and (9) show that the error estimates for solutions $\|\Delta y_\alpha\|/\|y\|$ are the same for the original and model examples under the condition of identity of $\|\widetilde{A}\|\eta$.

**An example from spectroscopy**

Let us illustrate the foregoing method of modeling (training, learning) by an example from the inverse spectroscopy problem (cf. [8]). The problem is to restore a spectrum via solving the Fredholm integral equation of the first kind (an ill-posed problem)

$$Ay \equiv \int_a^b K(\lambda,\lambda') y(\lambda') d\lambda' = f(\lambda), \quad c \leq \lambda \leq d, \quad (16)$$

where $K(\lambda,\lambda')$ is the spread function of a spectral device, $y(\lambda)$ is the true (desired) spectrum, $f(\lambda)$ is the measured (experimental) spectrum, $\lambda$ is the wavelength, $[a,b]$ are the limits for desired spectrum, $[c,d]$ are the limits for measured spectrum.

We assume that, instead of exact $f$ and $K$, we have $\widetilde{f}$ and $\widetilde{K}$ such that $\|\widetilde{f}-f\|\leq\delta$, $\|\widetilde{A}-A\|\leq\xi$. Equation (16) is solved by the Tikhonov regularization method according to (2), where $A^* = A^T$.

At first, we consider the *original example P* with known measured noisy spectrum $\widetilde{f}(\lambda)$ (Fig. 1) on a uniform grid of nodes $\lambda = \lambda_{\min}, \lambda_{\min}+h, \ldots, \lambda_{\max}$, where $\lambda_{\min} = c = 450$ nm, $\lambda_{\max} = d = 650$ nm, $h = \Delta\lambda = \mathrm{const} = 1$ nm is the discretization step, and $n = (\lambda_{\max}-\lambda_{\min})/h = 200$ is the number of discretization steps in $\lambda$.

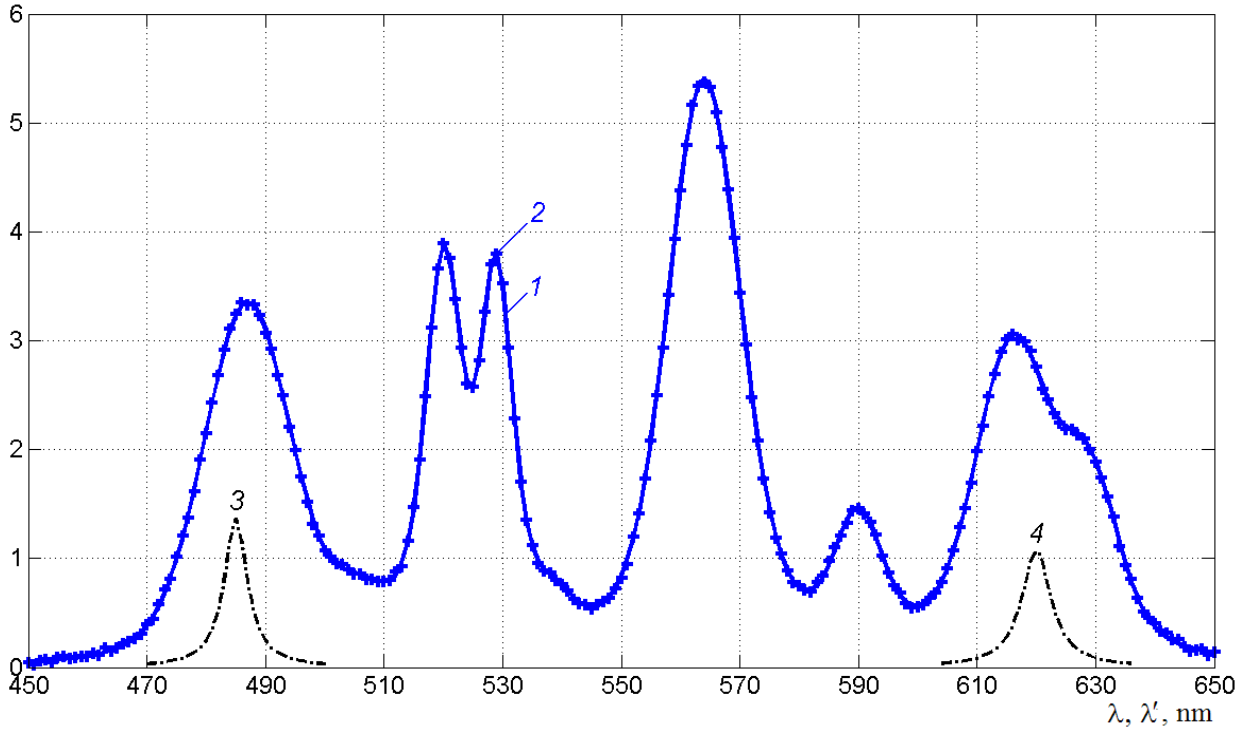

Fig. 1. Example *P*. *1* – exact spectrum $f(\lambda)$; *2* – noisy spectrum $\widetilde{f}(\lambda)$ and two cross-sections of SF: *3* – $10 K(485,\lambda')$ and *4* – $10 K(620,\lambda')$

It is assumed that the spread function (SF) $K = K(\lambda,\lambda')$ of spectral device has a variable width, i.e. is nondifference. As is known [8, 19], the SF width $w(\lambda)$ at level of 0,5 is proportional to wavelength $\lambda$. Therefore, we assume $w(\lambda) = q\lambda$, where $q = 0{,}015$. This corres-ponds to $w(c) = w(450 \text{ nm}) = 6{,}75 \text{ nm}$, $w(485 \text{ nm}) = 7{,}275 \text{ nm}$, $w(620 \text{ nm}) = 9{,}3 \text{ nm}$, and $w(d) = w(650 \text{ nm}) = 9{,}75 \text{ nm}$.

We use the dispersion SF

$$K(\lambda, \lambda') = \frac{w(\lambda)/2\pi}{(\lambda - \lambda')^2 + [w(\lambda)/2]^2}. \tag{17}$$

It is shown in [8] that SF of this type gives one of the most accurate restorations of a spectrum. To characterize the SF, along with the width $w(\lambda)$ at a level of 0,5, one may also use the integral SF width $W(\lambda)$ (the ratio of the SF area to its height)

$$W(\lambda) = \int_{-\infty}^{\infty} K(\lambda, \lambda') d\lambda' \Big/ K(\lambda, \lambda).$$

For a dispersion SF, we have: $W(\lambda) = (\pi/2)\, w(\lambda) \approx 1{,}571\, w(\lambda)$.

Figure 1 shows the SF $K(\lambda, \lambda')$ (17) at $\lambda = 485$ and 620 nm.

Analysis of Fig. 1 shows that the true (unknown) spectrum has, most likely, two close lines in the vicinity of $\lambda \approx 525$ nm and near $\lambda \approx 620$ nm, but they are poorly resolved in the measured spectrum $f(\lambda)$. Moreover, there is an indication that there is one weak line at $\lambda \approx 507$ nm, as well as at $\lambda \approx 543$ nm in true spectrum $y(\lambda)$. Thus, everything indicates that there are at least nine spectral lines in the spectrum $y(\lambda)$, although the number of lines in the measured spectrum $f(\lambda)$ is smaller (6 or 7).

In connection with this, the second (model, training) example $Q$ "close" to the original $P$ was modeled. The true spectrum of example $Q$ contains 9 as well as 8 and 10 spectral lines in the form of Gaussian (cf. [3, 8]), i.e. several examples $Q$ were modeled.

The measured spectra $f_Q(\lambda)$ in examples $Q$ were numerically calculated by formula

$$f_Q(\lambda) = \int_a^b K(\lambda, \lambda')\, y_Q(\lambda')\, d\lambda', \quad c \leq x \leq d.$$

Furthermore, $a = 460$ nm, $b = 640$ nm.

The measurement errors $\delta$ of the spectrum $f_P(\lambda)$ were estimated at about 1%, which corresponds to the standard deviation (SD) $\approx 0{,}02$. Therefore, the values of $f_Q(\lambda)$ were noisy by random errors with SD from 0,01 to 0,04, which corresponds to $\delta_{rel} \approx 0{,}5 - 2\%$ (because the value $\delta_{rel}$ is known inexactly). The SF in example $Q$ was taken in the form (17), moreover (since the SD is also known inexactly), $w(\lambda)$ was taken to be $w(\lambda) = q(1+\zeta)\lambda$, where $\zeta \in [-0{,}02, 0{,}04]$, which corresponds to $\xi_{rel} \approx 0 - 4\%$.

Further, the "close" model examples $Q$ were solved by the quadrature method with Tikhonov regularization via solving equation (2) at SF (17) for several values of the regularization parameter $\alpha$. It was calculated the dependence of the relative error of regularized solutions $y_\alpha(\lambda) = y_{\alpha Q}(\lambda)$ with respect to exact solutions $y(\lambda) = y_Q(\lambda)$:

$$\sigma_{rel}(\alpha) = \frac{\|y_\alpha(\lambda) - y(\lambda)\|}{\|y(\lambda)\|}.$$

Figure 2 shows dependences $\sigma_{rel}(\alpha)_Q$ for series of "close" model examples and for several values of errors $\delta_{rel} = \delta/\|f\|$ and $\xi_{rel} = \xi/\|A\|$ (the region between the curves *1* and *2*). Note that the curve *3* in the modeling method is supposed to be unknown; it is given for illustration only.

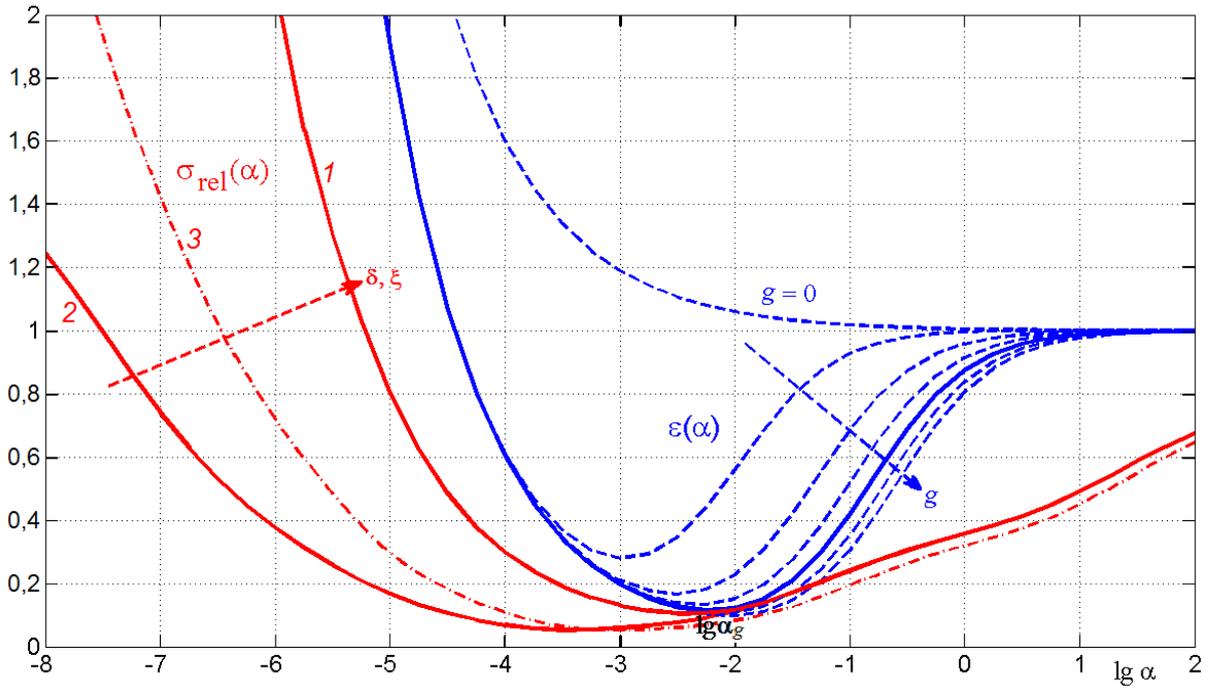

Fig. 2. Relative errors $\sigma_{rel}(\alpha)$ for examples $Q$ (*1* and *2* – the boundaries above and below, *3* – $\sigma_{rel}(\alpha)_P$) and envelope curves $\varepsilon(\alpha)$ for several values of parameter $g$

Figure 2 shows also several envelopes $\varepsilon(\alpha)$ according to (9) at $\|\widetilde{A}\|=\|A\|=0{,}843$ and $\eta=\delta_{rel}+\xi_{rel}=2\cdot10^{-2}$ for several values of $g$ from 0 to 0,1. We choose such value of $g$ at which one of curves $\varepsilon(\alpha)$ contacts set of curves $\sigma_{rel}(\alpha)$, namely, $g=0{,}045$. This corresponds to regularization parameter $\alpha=\alpha_g=10^{-2{,}2}$. It is seen from Fig. 2 that, despite the scatter of curves $\sigma_{rel}(\alpha)$ and $\varepsilon(\alpha)$, the value of $g$ and, as a consequence, $\alpha$ are estimated reliably.

Figure 3 shows solution (restored spectrum) at $\alpha_g=10^{-2{,}2}$, $\sigma_{rel}(\alpha_g)=0{,}073=7{,}3\%$. We can see that the spectrum is restored accurately: close lines are resolved and weak lines are separated.

**Remark 3.** Although we assume in the modeling method that the exact spectrum (solution) $y(\lambda)$ is unknown in original example *P*, we adduce the exact spectrum $y_P(\lambda)$ (Fig. 3) to demonstrate the potential possibilities of a technique. However, the spectrum $y_P(\lambda)$ is not used for choosing $\alpha_g$.

## Conclusion

There are a number of ways for choosing the regularization parameter $\alpha$ and estimating error $\|\Delta y_\alpha\|$ for regularized solution $y_\alpha$. Note the discrepancy principle [18], the generalized discrepancy principle [9], the modified discrepancy principle (the Raus–Gfrerer rule) [20], the cross-validation method [21], the L-curve criterion [22], the local regularizing algorithm [23], the new criterion of a posteriori choosing [15], the adaptive specialized generalized discrepancy principle [13], the new version of a posteriori choosing [16, 17], etc.

In these methods, as additional information about the solution, the errors $\delta$ and $\xi$ (as well as sourcewise representability of the solution) are usually used. As a result, the regularization parameter $\alpha$ at finite $\delta$ and $\xi$ is chosen reliably (but with some overstatement in comparison with $\alpha_{opt}$). Furthermore, the solution error $\|\Delta y_\alpha\|$ is obtained mainly in the form of asymptotic estimates, and the estimate $\|\Delta y_\alpha\|$ at finite $\delta$ and $\xi$ is usually obtained with a large overstatement (see. Fig. 2, the curve $g = 0$).

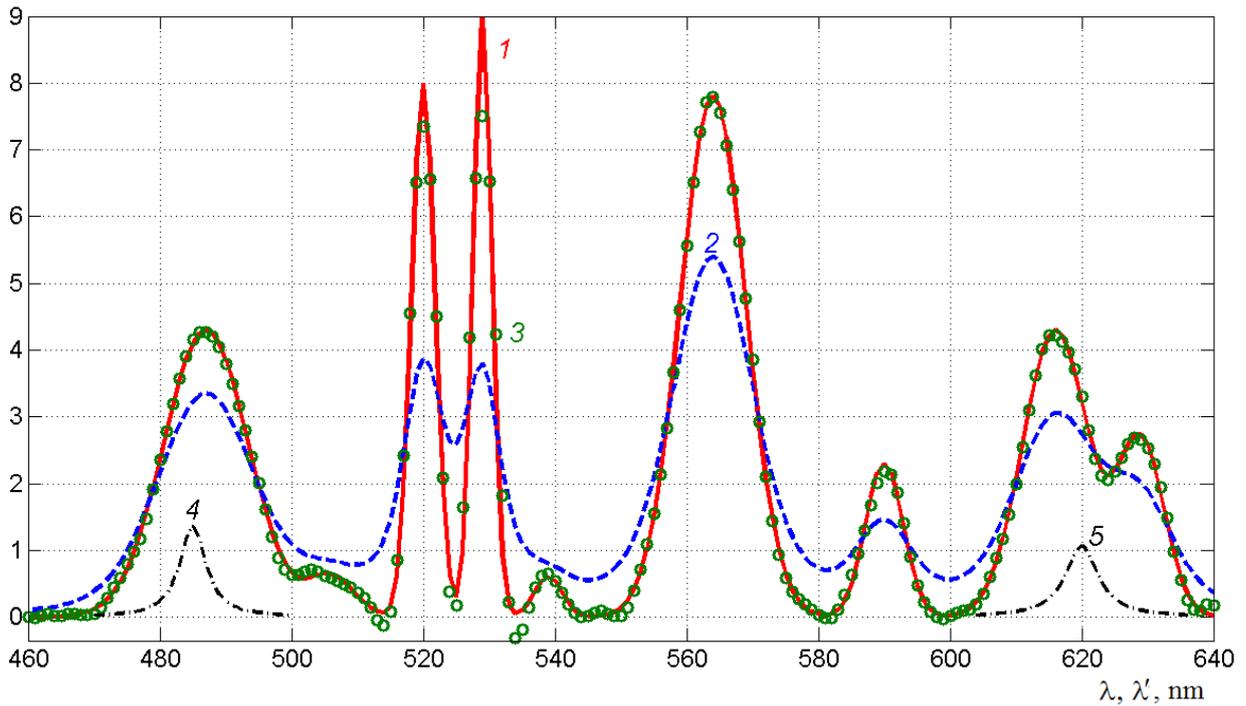

Fig. 3. Example *P*. *1* – true spectrum $y_P(\lambda)$; *2* – measured spectrum $f_P(\lambda)$;
*3* – restored spectrum $y_{\alpha P}(\lambda)$ at $\alpha = \alpha_g = 10^{-2.2}$ and two cross-sections of SF:
*4* – $10\,K(485,\lambda')$ and *5* – $10\,K(620,\lambda')$

In this paper, we develop a method of modeling or a method of training examples, which allows to choose $\alpha$ and, most importantly, to obtain unoverstated error estimate $\|\Delta y_\alpha\|$ (see Fig. 2, the contact of curves $\sigma_{rel}(\alpha)$ and $\varepsilon(\alpha)$ at $\alpha = \alpha_g$, as well as Fig. 3, the curve $y_\alpha(\lambda)_P$).

This work was supported by the Russian Foundation for Basic Research – RFBR (grant № 13-08-00442).

***Sizikov Valery Sergeevich*** – Saint-Petersburg National Research University of Information Technologies, Mechanics and Optics (ITMO University), Dr. Tech. Sci., Professor, sizikov2000@mail.ru

***Stepanov Andrey Valerievich*** – Vernadsky Crimean Federal University, Simferopol, Republic of Crimea, Dr. Tech. Sci., Professor, abc17101@yandex.ua